\documentclass[twocolumn]{autart}

\usepackage{times}
\usepackage{helvet}
\usepackage{courier}
\usepackage[T1]{fontenc}
\usepackage{amsmath}
\usepackage{amsfonts}
\usepackage{marvosym}
\usepackage{textcomp}
\usepackage{psfrag}
\usepackage{calc}
\usepackage{graphicx}
\usepackage[normalem]{ulem}
\usepackage{algorithm}
\usepackage{latexsym}

\newtheorem {prp} {Proposition}[section]
\newtheorem {dfn} {Definition}[section]

\def\IR{\mathbb{I{\kern -0.08em}R}}

\begin{document}

\begin{frontmatter}

\title{A sufficient condition to test identifiability of a nonlinear delayed-differential model with constant delays and multi-inputs}

%


\author{Carine Jauberthie}\ead{cjaubert@laas.fr},
\author{Louise Trav\'e-Massuy\`es}\ead{louise@laas.fr}

\address{CNRS; LAAS; 7, avenue du Colonel Roche, F-31077 Toulouse, France \\
Universit\'e de Toulouse; UPS, INSA, INP, ISAE; LAAS; F-31077 Toulouse, France\\}







\begin{abstract}
 In this paper, an original result in terms of a sufficient condition to test identifiability of nonlinear delayed-differential models with constant delays and multi-inputs is given. The identifiability is studied for the linearized system and a criterion for linear systems  with constant delays is provided, from which the identifiability of the original nonlinear system can be proved. This result is obtained by combining a classical identifiability result for nonlinear ordinary differential systems due to M.S. Grewal and K. Glover \cite{grewal} with the identifiability of linear delayed-differential models developed by Y. Orlov \emph{et al.} \cite{youri1}. 
This paper is a generalization of \cite{lilianne}, which deals with the specific case of nonlinear delayed-differential models with two delays and  a single input.
\end{abstract}

\end{frontmatter}

\section{Introduction}

Differential systems with delays enter into the modeling of many problems and are frequently used in domains like electronics, telecommunications, biology, epidemiology or aerospace. Satelite remote control or network communication protocol models fall for instance into this framework.
When the identifiability analysis of delays arises, it is 
impossible to directly use the classical criteria based on a similarity transformation approach \cite{vajda} or a series expansion approach \cite{walter} because delays occur in an implicit way in the state equations (argument of the functions of input and state).  
Then two approaches for identifiability analysis can be considered. \\
The first one consists in approximating the functions with delays \cite{carine} so that the approximate system is described by ordinary differential equations.  It is clear that the identifiability of this approximated system does not imply the identifiability of the original system.
The second approach consists in approximating the nonlinear system by linearization around an equilibrium state, which is the method followed in this paper. We extend a classical identifiability  result for nonlinear ordinary differential systems due to M.S. Grewal and K. Glover \cite{grewal} to nonlinear delayed-differential models.  We are hence interested in 
the identifiability of the linearized system and present a criterion to test the identifiability of a linear system  with constant delays \cite{lotfi}, \cite{youri1}, \cite{youri2}, from which the identifiability of the original nonlinear system can be proved. The result of \cite{grewal} is extended in two ways, by taking into account any number of delays acting on the state and on the input, and by considering the general case of multiple inputs, which was not considered in the extention proposed in \cite{lilianne}. The new proof requires to define a different norm for the inputs, to provide a bound for the estimation of the state norm that accounts for the delays and to deal with higher dimension.  \\
This paper is organized as follows.  In Section 2, we give some generalities on delayed-differential models and define identifiability for such systems. In Section 3, the result of identifiability with linearization around an equilibrium state is given through a criterion allowing us to test the identifiability of a linear system with constant delays \cite{lotfi}, \cite{youri1}, \cite{youri2}. Section 4 provides an illustrative example for which our identifiability analysis method is fully developped. The last section concludes the paper and discusses potential application domains which could benefit from our result.  






\section{Problem formulation}

In this work, we consider linear  and 
nonlinear systems with real positive delays. These systems are characterized by the length of their 
memory, i.e. the largest of their delays. The memory, positive, can be finite or infinite. It is supposed finite. 
For such a system, the state at one time point $t$ is defined on one 
interval $[t ', \ t]$, where $t'$ depends on the delays. If $\tau_{m}$ and $\nu_{m}$ 
are the greatest delays on the state and the input, respectively, the knowledge of $x(t)$ on 
$[t_{0}-\tau_{m}, \ t_{0}]$ and of $u(t)$ on $[t_{0}-\nu_{m}, \ t]$ are 
necessary and sufficient to determine $x(t)$ for all 
$t\geq t_{0}$. Thus, for a system with delays, the initial state (or 
initial function) must be given for all $t$ in $[t_{0}-\tau_{m}, \ 
t_{0}]$ where $t_{0}$ is the initial time of observation.\\
The nonlinear time-delayed system $\Gamma$, of specific state and input delays $\tau_{i}$ and $\nu_{j}$, can be written in the general form: 
\begin{equation}\label{modelenonlineairearetard}
\left\{\begin{array}{rl}
\dot{x}(t,P)=&f(x(t,P),x(t-\tau_{1},P),...,x(t-\tau_{l},P),\\
& u(t),u(t-\nu_{1}),...,,u(t-\nu_{r}), p_{1},...,p_{p}),\\
y(t,P)=&Cx(t,P),\ t>0,\\
x_{0}(s,P)=&\phi(s,P),\ -\tau_{m}\leq s\leq0,\\
u_{0}(s)=&\eta(s),\ -\nu_{m}\leq s\leq0.
\end{array}\right.
\end{equation}
where $x(t,P)\in\mathbb{R}^{n}$ and $y(t,P)\in\mathbb{R}^{m}$ denote the
state variables
and the measured outputs respectively. The input $u$ is
piecewise continuous with values in $\mathbb{R}^{k}$, and the initial function $x_{0}$ is continuous in
$s\in[-\tau_{m},0]$. The parameters $\tau_{1}$,...,$\tau_{l}$, $\nu_{1}$,...,$\nu_{r}$ represent the delays to be estimated. $\tau_{m}$ and $\nu_{m}$ are respectively the maximum delay of $\tau_{i}$ and $\nu_{j}$.  The vector of delays is included in $\mathcal{U}_{d}\subset \mathbb{R}^{l+r}$. The unknown constants $p_{1}$,..., $p_{p}$ belonging to $\mathcal{U}_{p}\subset \mathbb{R}^{p}$  and gathered in the vector $P_s$ must also be estimated. Hence the estimation concerns  the parameter vector $P=(p_{1},...,p_{p}, \tau_{1},...,\tau_{l}, \nu_{1},...,\nu_{r}) \in \mathcal{U}_{\mathcal {P}} \subset \mathbb{R}^{p+l+r}$. The function $f(.,.,.,.,p_{1},...,p_{p})$  is real and twice continuously  differentiable for every $P\in \mathcal{U}_{\mathcal {P}}$ on $M$ (a connected open subset of 
$\mathbb{R}^{n}\times...\times\mathbb{R}^{n}\times\mathbb{R}^{k}\times...\times\mathbb{R}^{k}$ such that $(x(t,P),x(t-\tau_{1},P),...,x(t-\tau_{l},P),u(t),u(t-\nu_{1}),...,,u(t-\nu_{r}))\in M$ for every $P\in \mathcal{U}_{\mathcal {P}}$ and every $t\in[0,\ T]$), where $T$ is the time duration. 
 $\phi(.,P)$ (respectively $\eta(.)$) is a continuous function, bounded on $[-\tau_{m},\ 0]$ in $\mathbb{R}^{n}$ (respectively on $[-\nu_{m},\ 0]$ in $\mathbb{R}^{k}$).
The output is a linear function of the state ($C$ is a matrix of appropriate dimensions).\\
The general theory of  systems with delays is developed in 
\cite{hale}, \cite{timedelaysystem}. Many works concern the analysis and the control of linear delayed models but there are much less works about identifiability analysis. The identifiability of these models has been analysed with restrictive conditions on the structure of the system \cite{nakagiri} or \cite{verduyn}.
The identifiability conditions of transfer functions are provided with sufficient nonsmooth inputs in \cite{youri1}, \cite{youri2}. 
But to our knowledge there is no general result for solving the identifiability problem for general nonlinear delayed-differential models with unknown constant delays such as (\ref{modelenonlineairearetard}). The proposed approach for solving the identifiability problem of systems given by 
(\ref{modelenonlineairearetard}) relies on the analysis of 
the identifiability of a form linearized around an equilibrium state. The definition of model identifiability considered in this paper is the following:
\begin{dfn}\label{dfniden}
  The model $\Gamma$ given by (\ref{modelenonlineairearetard}) is globally (resp. locally) identifiable at $P\in \mathcal{U}_{\mathcal {P}}$ if there exists a control $u$ such that, for any
  $\tilde{P}\in \mathcal{U}_{\mathcal {P}}$, the equality
  $P=\tilde{P}$ follows from $y(t,\tilde{P})=y(t,P)$
   $\forall t\in[0,T]$ (resp. if
there exists an open neighbourhood $W$ of $P$ such that $\Gamma$ is globally identifiable at $P$ with $\mathcal{U}_{\mathcal {P}}$
restricted to $W$).\\
\end{dfn}
In most models, there exist atypical points
in $\mathcal{U}_{\mathcal {P}}$  for which the model is unidentifiable. To account for these singularities, the previous definition can be generically extended. $\Gamma$ is
said to be (globally) structurally identifiable if it is (globally)
identifiable at all $P\in \mathcal{U}_{\mathcal {P}}$ except at a subset of points of
zero measure in $\mathcal{U}_{\mathcal {P}}$. \\
The following section formulates the identifiability problem and summarizes some results that are then used to derive our result.


\section{Identifiability results}

To assess the identifiability of system 
(\ref{modelenonlineairearetard}), our approach relies on testing the identifiability of a form  linearized around an equilibrium state, which leads to a sufficient condition for the nonlinear system.  The test is performed as proposed by MS. Grewal and K. Glover in \cite{grewal} for ordinary differential systems and it relies on the result of \cite{lotfi} in which L. Belkoura et al. propose a criterion allowing to test identifiability of a linear system with delays.\\
In the following section, we recall some definitions and results of identifiability for linear systems.

\subsection{Identifiability result for linear delayed-differential systems}
The identifiability of linear delayed-differential systems has been analyzed by  \cite{nakagiri}, \cite{verduyn} and \cite{lotfi}. The results of this subsection are mainly taken from \cite{lotfi}. The linear time-delayed model is supposed to be given by: 
\begin{equation}\label{lotfi1}
\left\{\begin{array}{rl}
\dot{x}(t)=&A_{0}x(t)+\sum_{i=1}^{l}A_{i}x(t-\tau_{i})\\
&+B_{0}u(t)+\sum_{i=1}^{r}B_{i}u(t-\nu_{i}),\\
x_{0}(s)=&\phi(s),\;-\tau_{l}\leq s\leq0,\\
u_{0}(s)=&\eta(s),\;-\nu_{r}\leq s\leq0 
\end{array}\right.
\end{equation}
with $0<\tau_{1}<\tau_{2}<...<\tau_{l}$ and $0<\nu_{1}<\nu_{2}<...<\nu_{r}$,
$x(t)\in\mathbb{R}^{n}$, $u(t)\in\mathbb{R}^{m}$ is piecewise continuous, $\phi\in\mathbb{L}^{2}(-\tau_{l},0;\mathbb{R}^{n})$. The delays belong to $\mathcal{U}_{d}$, a subset of $\mathbb{R}^{l}\times \mathbb{R}^{r}$. The components of the state are assumed to be available. Identifiability analysis is driven by the following definition:
\begin{dfn}
If the model is given by:
\begin{displaymath}
\left\{\begin{array}{rl}
\dot{\tilde{x}}(t)=&\tilde{A}_{0}\tilde{x}(t)+\displaystyle \sum_{i=1}^{l}\tilde{A}_{i}\tilde{x}(t-\tilde{\tau}_{i})+\tilde{B}_{0}u(t)\\
& +\displaystyle \sum_{i=1}^{r}\tilde{B}_{i}u(t-\tilde{\nu}_{i}),\\
\tilde{x}_{0}(s)=&\tilde{\phi}(s),\;-\tilde{\tau}_{l}\leq s\leq0,\\
u_{0}(s)=&\eta(s),\;-\tilde{\nu}_{r}\leq s\leq0,
\end{array}\right.
\end{displaymath}
then the matrices of coefficients $A_{i},\;(i=0,...,l)$, $B_{j},\;(j=0,...,r)$ and the delays $\tau_{k}$, $(k=1,...,l)$, $\nu_{p},\;(p=1,...,r)$ 
of system (\ref{lotfi1}) are structurally  globally identifiable under $u(t)$ if:
$x(t)=\tilde{x}(t),t\geq0\;\Rightarrow\left\{\begin{array}{l}
A_{i}=\tilde{A}_{i},\; (i=0,...,l),\\
\tau_{k}=\tilde{\tau}_{k}\;(k=1,...,l),\\
B_{j}=\tilde{B}_{j},\;(j=0,...,r),\\
\nu_{p}=\tilde{\nu}_{p},\;(p=1,...,r).
\end{array}\right.
$
\end{dfn}
In the following, for $z\in\mathbb{C}$, we note $A(z)$ and $B(z)$ the
expressions $A(z)=A_{0}+A_{1}z^{\tau_{1}}+...+A_{l}z^{\tau_{l}},$ and 
$B(z)=B_{0}+B_{1}z^{\nu_{1}}+...+B_{r}z^{\nu_{r}}.$
Then we have the following result.
\begin{thm}\label{thmlotfi}
(From  \cite{lotfi}) Considering the model (\ref{lotfi1}), if there exists a complex number $z$ such that:
$$rank[B(z)|A(z)B(z)|...|A^{n-1}(z)B(z)]=n,$$
then there exists a control $u(t)$ for which the matrices of
coefficients $A_{i},(i=0,...,l)$, $B_{j}$, $(j=0,...,r)$ and the delays
$\tau_{k},(k=1,...,l)$ et $\nu_{p},(p=1,...,r)$ are structurally globally identifiable.
\end{thm}
In the case of multi-inputs, it is sufficient to find an input of type `` square pulse '' with discontinuities of full rank, i.e. discontinuities are incommensurable and verify a rank condition \cite{lotfi}. Otherwise, a piecewise constant input can be appropriate.\\
The above result is used in Section 3.2 to test the identifiability of the linearized nonlinear delayed-differential system.

\subsection{A sufficient condition for identifiability of nonlinear
delayed-differential systems}

The identifiability of nonlinear dynamical systems (without delays) by linearization around an operating state 
has been shown by M.S. Grewal and K. Glover \cite{grewal}. The idea 
of the present paper is to
extend this approach to the class of nonlinear delayed-differential models with unknown constant delays. The considered operating state of system (\ref{modelenonlineairearetard}) is an 
equilibrium state $x_{e}$ corresponding to a constant input $\bar u$:
\begin{equation}
\label{equ}
\left\{ \begin{array}{l}
0=f(x_{e},...,x_{e},\bar u,...,\bar u,P_s)\,,\\
y_{e}=Cx_{e}\,.
\end{array}
\right.
\end{equation}
To simplify the notation, the point $(x_{e},...,x_{e},\bar u,...,\bar u,P_s)$ given in (\ref{equ}) is noted $E$.\\ Given a constant input $\bar u\in\mathbb{R}^{k}$, let us consider the
operating state $x_{e}$ defined by (\ref{equ}). For a given $P\in{\mathcal
U}_P$, the state $x_{e}$ may not be unique.

Moreover $x(t,P)=x_{e}$ is the solution of:
\begin{equation}
\label{eq1}
\left\{ \begin{array}{l}
\dot{x}(t,P)=f(x(t,P),x(t-\tau_{1},P),...,x(t-\tau_{l},P),\\
\hspace{1cm} u(t),u(t-\nu_{1}),...,u(t-\nu_{r}),P_s)\,,\;
t\in[0,T]\,,\\
x_{0}(s,P)=x_{e}\,,\;s\in[-\tau_{m},0]\,,\\
y(t,P)=C x(t,P)\,.
\end{array}
\right.
\end{equation}
corresponding to the input $u(t)=\bar u$, $t\in[-\nu_{m},T]$.\\
Now let us consider a perturbation $\nu(t)$ of $\bar u$ for $t\geq 0$
and let us define the corresponding system:
\begin{equation}
\label{eq2}
\left\{ \begin{array}{rl}
\dot{x}(t,P)=&f(x(t,P),x(t-\tau_{1},P),...,x(t-\tau_{l},P),\\
 & \bar u+\nu(t),\bar u+\nu(t-\nu_{1}),...,\\
& \bar u+\nu(t-\nu_{r}),P_s)\,,\; t\in[0,T]\,,\\
x_{0}(s,P)=&x_{e}\,,\;s\in[-\tau_{m},0]\,,\\
y(t,P)=&C x(t,P)\,.
\end{array}
\right.
\end{equation}
If we introduce $x_{\delta} (t,P)=x(t,P)-x_{e}$, the
previous system can be expressed as:
\begin{equation}
\label{eq3}
\left\{ \begin{array}{rl}
\dot{x}_{\delta}(t,P)=&f(x_{e}+x_{\delta}(t,P),x_{e}+x_{\delta}(t-\tau_{1},P),...,\\
& x_{e}+x_{\delta}(t-\tau_{l},P),\bar u+\nu(t), \bar u+\nu(t-\nu_{1}),...,\\
& \bar u+\nu(t-\nu_{r}),P_s)\,,\;
t\in[0,T]\,,\\
x_{\delta}(s,P)=&0\,,\;s\in[-\tau_{m},0]\,,\\
y_{\delta}(t,P)=&C x_{\delta}(t,P)\,.
\end{array}
\right.
\end{equation}
If the variables of $f$ introduced in (\ref{eq3}) are denoted $z_{i}=x_{e}+x_{\delta}(t-\tau_{i},P)$, $i=0,...,l$ where $\tau_{0}=0$ and 
$w_{j}=\bar u+\nu(t-\nu_{j})$,  $j=0,...,r$ where $\nu_{0}=0$, we get: 
$\begin{array}{rrl} 
f: &
\mathbb{R}^{n}\times...\times\mathbb{R}^{n}\times\mathbb{R}^{k}\times...\times\mathbb{R}^{k}\times{\mathcal U}_p &
\to \mathbb{R}^{n} \\
  & (z_{0},..., z_{l},w_{0},..., w_{r}, P_{s}) &\mapsto f(z_{0},..., z_{l},\\
& & \hspace{0.5cm} w_{0},...,w_{r}, P_{s})
   \end{array}
$
and if the following matrices are introduced:
\begin{equation}
     \left\{\begin{array}{l}
     A_{i}(P)=\nabla_{z_{i}}f(E),\ \ i=0,...,l,\\
    
     B_{j}(P)=\nabla_{w_{j}}f(E),\ \ i=0,...,r,\\
   
     \end{array}\right.
\end{equation}
where $\nabla_{z_{i}}f(E)$ ($\nabla_{w_{i}}f(E)$) represents the jacobien matrix of $f$ with respect to $z_{i}$ ($w_{i}$) calculated at $E=(x_{e},...,x_{e},\bar u,...,\bar u,P_{s})$, then the linear delayed system, obtained by the linearization of
(\ref{modelenonlineairearetard}), is given by:
\begin{equation}
\label{eq4}
\left\{ \begin{array}{l}
\dot{\xi}(t,P)=\Sigma_{i=0}^{l}A_{i}(P)\xi(t-\tau_{i},P)+\Sigma_{j=0}^{r}B_{j}(P)\nu(t-\tau_{j})\,,\\
\xi(s,P)=0\,,\;s\in[-\tau_{m},0]\,,\\
\eta(t,P)=C \xi(t,P)\,.
\end{array}
\right.
\end{equation}
This last system is of the form (\ref{lotfi1}) when $C=\mathbb{I}$, and its identifiability is obtained by Theorem 1. \\
The following proposition is an extension of the results of
\cite{grewal} to delayed-differential systems.
\begin{prp}
     \label{prop}
If the model (\ref{eq4}) is (structurally) globally (resp. locally) identifiable
at $P\in{\mathcal U}_P$ (resp. $P\in W$)\footnote{Let us recall that $W$ is an open neighbourhood of $P$.}, then the model 
(\ref{eq1}) is  (structurally) globally (resp. locally) identifiable
at $P$, with a perturbation of $\bar{u}$ as input.
\end{prp}
The principle of the proof, as developped below, is that if $\bar{\nu}$ is a control providing the identifiability of (\ref{eq4}), then the identifiability of the model (\ref{eq1}) can be obtained by the control $\bar{u}+\varepsilon \bar{\nu}/|| \bar{\nu} ||$, with $\varepsilon$ small. This proof requires to define an appropriate norm for the inputs, to provide a bound for the estimation of the state norm that accounts for the delays.  The proof is given in the general case in which $C$ is not necessarily equal to $\mathbb{I}$. However, if this is the case, a more general result than Theorem \ref{thmlotfi} is required to assess the identifiability of the linearized system. 

{\it Proof} - Let $\tilde P\neq P$, $\tilde P\in \mathcal{U}_P$ (resp. $\tilde P\in W$). Let us consider the equilibrium states $x_{e}$ and $\tilde x_{e}$ corresponding respectively to $\tilde P$ and $P$:
\begin{equation}
\label{eq5}
\left\{ \begin{array}{ll}
0 & =f(E)\,,\\
y_{e} & =Cx_{e}\,,
\end{array}
\right.
\left\{ \begin{array}{ll}
0 & =f(\tilde E )\,,\\
\tilde y_{e} & =C\tilde x_{e}\,.
\end{array}
\right.
\end{equation}
Two cases are possible:

\hspace{0.5cm} - $y_{e}\neq \tilde y_{e}$, then the system
(\ref{eq1}) is globally (resp. locally) identifiable
at $P$ since there exists an input $\bar u$ such that the outputs
of the system (\ref{eq1}) $y_{e}$ and $\tilde y_{e}$ are distinct.

\hspace{0.5cm} - $y_{e}=\tilde y_{e}$. Let us evaluate the difference
between the outputs of the system (\ref{eq1}):
$$\begin{array}{rl}
y(t,P)-y(t,\tilde P)  =&
y(t,P)-y_{e}-(y(t,\tilde P)-\tilde y_{e})\\
   =&C(x(t,P)-x_{e})-C(x(t,\tilde P)-\tilde x_{e})\\
    =&C(x_{\delta}(t,P)-x_{\delta}(t,\tilde P))\\
    =&C(x_{\delta}(t,P)-\xi(t,P))+\\
& 
   C(\xi(t,P)-\xi(t,\tilde P))-\\
& C(x_{\delta}(t,\tilde P)-
   \xi(t,\tilde P)).
\end{array}$$
Hence
\begin{equation}\label{eq6}
\begin{array}{rl}
   
   y(t,P)-y(t,\tilde P)=&\eta(t,P)-\eta(t,\tilde P)+ C((x_{\delta}(t,P)\\
 & -\xi(t,P)) - (x_{\delta}(t,\tilde P)-\xi(t,\tilde P))).
\end{array}
\end{equation}
Since the system (\ref{eq4}) is globally (resp. locally) identifiable
at $P$, there exists an input $\bar\nu(t)$ such that the
corresponding outputs $\bar\eta(t,P)$ and $\bar\eta(t,\tilde P)$
are distinct. Now, take the input:
\begin{equation}
     \label{nu}
     \nu(t)=\frac{\varepsilon}{\vert\vert \bar\nu\vert\vert }\bar\nu(t)
     \quad (\vert\vert \bar\nu\vert\vert=\vert\vert
     \bar\nu\vert\vert_{L^{2}(0,T)}),
\end{equation}
in which $\varepsilon$ is chosen judiciously in the following.\\
Since the system (\ref{eq4}) is linear, the output corresponding to
the input $\nu(t)$ is given by: 
$$\eta(t,P)=\frac{\varepsilon}{\vert\vert \bar\nu\vert\vert }
\bar\eta(t,P).$$
Consequently, the $L_{2}$-norm of the difference between the outputs
of the system (\ref{eq4}) is
\begin{equation}
     \label{eq7}
\vert\vert \eta(.,P)-\eta(.,\tilde P)\vert\vert=
    \frac{\varepsilon}{\vert\vert \bar\nu\vert\vert}
\vert\vert \bar\eta(.,P)-\bar\eta(.,\tilde P)\vert\vert
=\varepsilon K,
\end{equation}
where the constant $K$ is strictly positive.\\
The $L_{2}$-norm of (\ref{eq6}) leads to:
$$\begin{array}{rl}
\vert\vert y(.,P)-y(.,\tilde P)\vert\vert \geq & 
\vert\vert \eta(.,P)-\eta(.,\tilde P)\vert\vert \\
& -\vert\vert C((x_{\delta}(.,P)-\xi(.,P))\\
 & -(x_{\delta}(.,\tilde P)-\xi(.,\tilde P)))\vert\vert
\end{array}
$$
which, by introducing (\ref{eq7}), gives
\begin{equation}
     \label{eq7bis}
\begin{array}{rl}
\vert\vert y(.,P)-y(.,\tilde P)\vert\vert \geq & \varepsilon K
-\vert\vert C\vert\vert(\vert\vert 
x_{\delta}(.,P)-\xi(.,P)\vert\vert \\
& + \vert\vert
x_{\delta}(.,\tilde P)-\xi(.,\tilde P)\vert\vert).
\end{array}
\end{equation}
In the following of the proof, a lower bound is provided for $\vert\vert y(.,P)-y(.,\tilde P)\vert\vert$, which requires to estimate $\vert\vert  x_{\delta}(t,P)-\xi(t,P)\vert\vert_{\mathbb{R}^{n}}$.

\textit{Estimation of }$\vert\vert
     x_{\delta}(t,P)\vert\vert_{\mathbb{R}^{n}}$ (needed to estimate $\vert\vert
     x_{\delta}(t,P)-\xi(t,P)\vert\vert_{\mathbb{R}^{n}}$). The integration of (\ref{eq3}) leads to
\begin{equation}
     \label{eq8}
\begin{array}{rl}
x_{\delta}(t,P)=&\int_{0}^t\big[ 
f(x_{e}+x_{\delta}(s,P),x_{e}+x_{\delta}(s-\tau_{1},P),...,\\
& x_{e}+x_{\delta}(s-\tau_{l},P), \bar u+\nu(s),\bar u+\nu(s-\nu_{1}),...,\\
& \bar u+\nu(s-\nu_{r}),P_{s}-f(E)\big] ds\,.
\end{array}
\end{equation}
In the following, the terms $L^{(1)}_{i},\ L^{(2)}_{i}$ and $L_{i}$ denote real constants. The assumption of smoothness of the function $f$ gives:
\begin{equation}
     \label{eq9}
\begin{array}{rl}
\vert\vert x_{\delta}(t,P)\vert\vert_{\mathbb{R}^{n}}\leq &
     \int_{0}^t (
\Sigma_{i=0}^{l}L^{(1)}_{i}\vert\vert x_{\delta}(s-\tau_{i},P)\vert\vert_{\mathbb{R}^{n}}\\
&+\Sigma_{j=0}^{r}L^{(2)}_{j}\vert\vert\nu(s-\nu_{j})\vert\vert_{\mathbb{R}^{k}})ds
\end{array}
\end{equation}
\begin{equation}
     \label{eq10}
\vert\vert x_{\delta}(t,P)\vert\vert_{\mathbb{R}^{n}}\leq
     \int_{0}^t (L_{1}\vert\vert
    x_{\delta}(s,P)\vert\vert_{\mathbb{R}^{n}}+L_{2}\vert\vert\nu(s)\vert\vert_{\mathbb{R}^{k}})ds.
\end{equation}
This last estimation is due to the fact that when $s \in [0, t]$, $s-\tau_i \in [-\tau_i, t-\tau_i]$. But $x_{\delta}(s,P)=0$ for
$s\in[-\tau_{m},0]$. Then considering the values of $\vert\vert x_{\delta}(t,P)\vert\vert_{\mathbb{R}^{n}}$ for $s \in [0, t]$ and taking $L_1$ as an upper bound of the coefficients  $L^{(1)}_{i}, i= 0, \dots, l$, provides an upper bound of the first term of the integral of (\ref{eq9}). The same reasoning stands for the second term, given that $\nu(s)=0$ for $s\in[-\nu_{m},0]$.
Moreover $\vert\vert\nu\vert\vert=\varepsilon$, thus by applying Gronwall's lemma it leads to:
\begin{equation}
     \label{eq11}
\vert\vert x_{\delta}(t,P)\vert\vert_{\mathbb{R}^{n}}\leq
    L_{3}\varepsilon\,.
\end{equation}
\textit{Estimation of }$\vert\vert
     x_{\delta}(t,P)-\xi(t,P)\vert\vert_{\mathbb{R}^{n}}$. Let us recall that:
\begin{equation}
     \label{eq12}
     \begin{array}{l}
\displaystyle\frac{d}{dt}\left(x_{\delta}(t,P)-\xi(t,P)\right)= \\
f(x_{e}+x_{\delta}(t,P),x_{e}+x_{\delta}(t-\tau_{1},P),...,x_{e}+x_{\delta}(t-\tau_{l},P),\\
 \bar u+\nu(t),\bar u+\nu(t-\nu_{1}),...,\bar u+\nu(t-\nu_{r}),P_{1},...,P_{p})\\
-(\Sigma_{i=0}^{l}A_{i}(P)\xi(t-\tau_{i},P)+\Sigma_{j=0}^{r}B_{j}(P)\nu(t-\tau_{j})).
\end{array}
\end{equation}
The assumption of smoothness of the function $f$ gives:
\begin{multline}
     \label{eq13}
\frac{d}{dt}\left(x_{\delta}(t,P)-\xi(t,P)\right)=
\Sigma_{i=0}^{l}A_{i}(P)(x_{\delta}(t-\tau_{i},P)\\
 -\xi(t-\tau_{i},P))+\mathcal{E}(t),
\end{multline}
where   \vspace{-0.9cm}  
\begin{equation}
\label{eq14a}
\begin{array}{rl}
\vert\vert \mathcal E(t)\vert\vert_{\mathbb{R}^{n}}\leq & L_{4} (
\Sigma_{i=0}^{l}\vert\vert x_{\delta}(t-\tau_{i},P)\vert\vert_{\mathbb{R}^{n}}^{2}\\
&+\Sigma_{j=0}^{r}\vert\vert \nu(t-\tau_{j})\vert\vert_{\mathbb{R}^{k}}^{2})
\end{array}
\end{equation}
and, by applying (\ref{eq11}): $\vert\vert \mathcal E(t)\vert\vert_{\mathbb{R}^{n}}\leq L_{5} \varepsilon^{2}$.

The same approach as in the previous estimation yields:
\begin{multline}
     \label{eq15}
\vert\vert x_{\delta}(t,P)-\xi(t,P)\vert\vert_{\mathbb{R}^{n}}\leq 
L_{6} \int_{0}^t \vert\vert 
x_{\delta}(s,P)-\xi(s,P)\vert\vert_{\mathbb{R}^{n}}ds \\
 +L_{7}\varepsilon^{2}
\end{multline}
    then (Gronwall's lemma)
$\vert\vert x_{\delta}(t,P)-\xi(t,P)\vert\vert_{\mathbb{R}^{n}}\leq
L_{8}\varepsilon^{2}.$

\textit{Underestimation of }$\vert\vert
     y(.,P)-y(.,\tilde P)\vert\vert$.\\ The overestimation of $\vert\vert x_{\delta}(t,P)-\xi(t,P)\vert\vert_{\mathbb{R}^{n}}$ applied to (\ref{eq7bis}) leads to 
$\vert\vert y(.,P)-y(.,\tilde P) \vert\vert\geq \varepsilon
K -L\varepsilon^{2}$.
Now, there exists $\varepsilon>0$ such that $\varepsilon
K -L\varepsilon^{2}>0$. Therefore, the outputs $y(.,P)$ and
$y(.,\tilde P)$ are distinct and the system (\ref{eq1}) is structurally  
globally (resp. locally) identifiable
at $P$.
\hfill$\square$

\section{Illustrative example}

In this section, the following example taken from \cite{ex} is used to illustrate the approach.
\begin{equation}\label{eqex}
\left\{\begin{array}{rl}
\dot{x}(t)=&-x(t)+(1+\sin^{2}(x(t)))y(t)+x^{2}(t-\tau_{1}),\\
\dot{y}(t)=&x(t)y(t)+v(t)+w(t)\\
& +(1+\sin^{2}(x(t)))u_{1}(t)+y(t-\tau_{2}),\\
\dot{v}(t)=&-v(t)+w(t)+v(t-\tau_{3}),\\
\dot{w}(t)=&(y(t)+v(t))w(t)-x(t)u_{1}(t) +(2-\sin(v(t)w(t)\\
& -x(t)))u_{2}(t) +x(t-\tau_{1})w(t-\tau_{4}).
\end{array}\right.
\end{equation}
 The identifiability result provided in Section 3.2 is fully developped and allows us to conclude on the identifiability of the system. In a second step, the system is parameterized by introducing thirteen constant parameters and the identifiability analysis is performed similarly, allowing us to conclude to the identifiability of the parameterized system as well. 
We assume that the state variables are available. Let us consider inputs $u_{1_e}$ and $u_{2_e}$. An equilibrium state 
$(x_e,y_e,v_e,w_e)$ is given by the resolution of (\ref{eqex}) in which the derivatives are set to $0$ and the state variables and any of their delayed representatives are set to the corresponding equilibrium values $x_e,y_e,v_e,w_e$. 
The third equation then provides $w_e=0$  and studying the fourth equation on an appropriately chosen domain leads to the existence of  $x_e$ such that $-x_eu_{1_e}+(2+\sin(x_e))u_{2_e}=0$. Given $x_e$, we obtain $y_e=\displaystyle\frac{x_e-x_e^{2}}{1+\sin^{2}(x_e)}$ from the first equation. Finally, using $x_e$ and $y_e$ in the second equation, we get $v_e=-x_ey_e-(1+\sin^{2}(x_e))u_{1_e}-y_e$. \\
For a matrix $M$, let us note $M_{i,j}$ its components where $i$ and $j$ are the row and column numbers, respectively. After linearization of the model (\ref{eqex}), we obtain the matrices $A_{i}$, $i=0,1,2,3,4$, and $B_{0}$ of the system put in form (2):\\
\vspace{0.1cm}
\newline
\scalebox{0.9}{
$A_{0}=\left (\begin{array}{cccc}
-1+\sin(2x_{e})y_{e} & 1+\sin^{2}(x_e) & 0 & 0\\
y_e+\sin(2x_e)u_{1_e} & x_e &  1 & 1 \\
 0 & 0 & -1 & 1\\
u_{2_e}\cos(x_{e})-u_{1_e} & 0 & 0 & y_{e}+v_{e}p_{e}
\end{array}\right),$}\\
with $p_{e}=1-u_{2_e}\cos(x_{e})$ and,\\
\scalebox{0.9}{
$
\left\{\begin{array}{l}
A_{1_{i,j}}=2x_{e}\ {\rm{if}}\ (i,j)=(1,1),\ 0\ {\rm{otherwise}},\\
A_{2_{i,j}}=1\ {\rm{if}}\ (i,j)=(2,2),\ 0\ {\rm{otherwise}},\\
A_{3_{i,j}}=1\ {\rm{if}}\ (i,j)=(3,3),\ 0\ {\rm{otherwise}},\\
A_{4_{i,j}}=x_e\ {\rm{if}}\ (i,j)=(4,4),\ 0\ {\rm{otherwise}},\\
B_{0_{i,j}}=\left\{\begin{array}{l} 
1+\sin^{2}(x_e)\ {\rm{if}}\ (i,j)=(2,1),\ -x_{e}  \ {\rm{if}}\ (i,j)=(4,1),\\
2+\sin(x_{e})  \ {\rm{if}}\ (i,j)=(4,2),\ 0\ \rm{otherwise}.
\end{array}\right.
\end{array}\right.
$}\\
\newline
By using the symbolic toolbox of Matlab, one can check that the condition $rank[B(z)|A(z)B(z)|...|A^{3}(z)B(z)]=4$ is satisfied for several values of $z$, for instance $z=2$ or $z=1+i$. 
Then there exists a control $u(t)$ for which the delays
$\tau_{k}, k=1,2,3,4,$ are structurally globally identifiable. Thus, by Proposition 3.1, the original nonlinear system (\ref{eqex}) is structurally globally identifiable.\\
Let us now consider the parameterized system: 
\begin{equation}\label{eq_ex2}
\left\{\begin{array}{rl}
\dot{x}(t)=&-p_{1}x(t)+(1+p_{2}\sin^{2}(x(t)))y(t)+p_{3}x^{2}(t-\tau_{1}),\\
\dot{y}(t)=&p_{4}x(t)y(t)+p_{5}v(t)+p_{6}w(t)\\
& +p_{7}(1+\sin^{2}(x(t)))u_{1}(t)+p_{8}y(t-\tau_{2}),\\
\dot{v}(t)=&p_{9}(v(t-\tau_{3})-v(t))+p_{10}w(t) ,\\
\dot{w}(t)=&p_{11}(y(t)+v(t))w(t)+p_{12}[(2-\sin(v(t)w(t)\\
& -x(t)))u_{2}(t)-x(t)u_{1}(t)]+p_{13}x(t-\tau_{1})w(t-\tau_{4}).
\end{array}\right.
\end{equation}
This system is obtained from system (\ref{eqex}) by introducing thirteen constant parameters $p_{i}$, $i=1,...,13$. These parameters are constants to be identified, given in the vector $P_s$. Assuming that $p_{10}$ and $p_{12}$ are different from 0, the same reasoning as for system (\ref{eqex}) leads to an equilibrium state $(x_e,y_e,v_e,w_e)$, in which $w_e=0$, $y_e$ and $v_e$ depend on the parameters, but not $x_e$. After linearization around this equilibrium state, we obtain the following matrices:\\
\newline
\scalebox{0.9}{
$A_{0}=\left (\begin{array}{cccc}
A_{0_{11}} & A_{0_{12}} & 0 & 0 \\
A_{0_{21}} & x_{e}p_{4} &  p_{5} &  p_{6}\\
0 & 0 & -p_{9} & p_{10} \\
A_{0_{41}} & 0 & 0 & A_{0_{44}}
\end{array}\right),
\left \{\begin{array}{l}
A_{0_{11}}= -p_{1}+\sin(2x_{e})y_{e}p_{2},\\ 
A_{0_{12}}= 1+\sin^{2}(x_e)p_{2},\\
A_{0_{21}}= y_{e}p_{4}+\sin(2x_e)u_{1_e}p_{7},\\
A_{0_{41}}= p_{12}(h_{e}-u_{1_e}),\\
A_{0_{44}}= p_{11}(y_{e}+v_{e})-v_{e}h_{e}p_{12},
\end{array}\right.$}
with $h_{e}=u_{2_e}\cos(x_e)$ and,

\scalebox{0.9}{
$
\left\{\begin{array}{l}
A_{1_{i,j}}=2x_{e}p_{3}\ {\rm{if}}\ (i,j)=(1,1),\ {\rm{and}}\ 0\ \rm{otherwise},\\
A_{2_{i,j}}=p_{8}\ {\rm{if}}\ (i,j)=(2,2)\ {\rm{and}}\ 0\ {\rm{otherwise}},\\
A_{3_{i,j}}=p_{9}\ {\rm{if}}\ (i,j)=(3,3)\ {\rm{and}}\ 0\ {\rm{otherwise}},\\
A_{4_{i,j}}=x_{e}p_{13}\ {\rm{if}}\ (i,j)=(4,4)\ {\rm{and}}\ 0\ {\rm{otherwise}},\\
B_{0_{i,j}}=\left\{\begin{array}{l}
(1+\sin^{2}(x_e))p_{7}\ {\rm{if}}\ (i,j)=(2,1),\\
-x_{e}p_{12}\ {\rm{if}}\ (i,j)=(4,1),\\
 (2+\sin(x_{e}))p_{12}\ {\rm{if}}\ (i,j)=(4,2)\ {\rm{and}}\ 0\ {\rm{otherwise}}.
\end{array}\right.
\end{array}\right.
$}

%
In the case of system (\ref{eq_ex2}), after linearization, the matrices $A_{i}$, $i=0,...,4$ and $B_{0}$ depend on the vector of parameters, i.e. $A_{i}=A_{i}(P_s)$ and $B_{0}=B_{0}(P_s)$. Hence the identifiability result is not immediate. We first prove the identifiability of the matrices and delays with the rank condition. We obtain the existence of a control such that the matrices of coefficients and delays are structurally globally identifiable. We then check whether this implies $P=\tilde{P}$. More precisely, we obtain $A_{i}(P)=A_{i}(\tilde{P})$, $i=0,...,4$ and $B_{0}(P)=B_{0}(\tilde{P})$ which leads, for example for $A_{1}(P)_{11}=A_{1}(\tilde{P})_{11}$, to $2x_{e}p_{3}=2x_{e}\tilde{p}_{3}$. As $x_{e}$ is known and different from 0, we obtain the identifiability of $p_{3}$. Repeating the same reasoning, we obtain the structural global  identifiability of all parameters. Consequently, the original nonlinear system (\ref{eq_ex2}) is structurally globally  identifiable.

\section{Conclusion}

This paper proves a sufficient condition for the identifiability of nonlinear
delayed-differential systems by linearization around an operating point. It generalizes the work of \cite{lilianne}, which deals with the specific case study of a nonlinear delayed-differential model with two delays and a single input. The proposed condition has been tested on two relevant examples corresponding to the considered type of model. \\
The condition that has been exhibited is only sufficient and further work should concentrate on finding the necessary part. It would be particularly interesting to find necessary and
sufficient conditions for the identifiability of nonlinear delayed
systems directly from the input-output relations with
parameters like in the case of nonlinear ordinary differential systems. \\
Identifiability is an important property that determines the system-based approach of control theory in which most of the modeling is performed by estimating the parameters of an priori given model structure. The domains in which non linear phenomena need to be represented are numerous and call for the kind of models that are considered in this paper. This is the case in the aerospace domain for which we are particularly interested in fault tolerant control. In this context, identifiability is a condition for on-line applications. It is indeed critical to detect the fault and immediatly identify the corresponding fault model, so that the control laws can be reconfigured appropriately. This exemplifies a research field that can certainly benefit from the results presented in this paper.


\end{document}